\documentclass[12pt]{article}
\usepackage{amsmath,amssymb,amsfonts,graphicx,wrapfig}

\newtheorem{theorem}{Theorem}
\newtheorem{lemma}{Lemma}
\newtheorem{corol}{Corollary}

\newcommand{\Pois}{{\mathrm{Pois}}}
\newcommand{\FO}{{\mathrm{FO}}}
\newcommand{\MSO}{{\mathrm{MSO}}}
\newcommand{\tr}{{\mathrm{trees}}}
\newcommand{\gr}{{\mathrm{graphs}}}
\newcommand{\EHR}{{\mathrm{EHR}}}

\begin{document}

\begin{center}
{\Large
Monadic second-order properties of very sparse random graphs~\footnote{M.E. Zhukovskii is supported by the grant 16-11-10014 of Russian Science Foundation. He proposed ideas of the proofs of the main results. Moreover, he proved Theorem 6.}\\
}
\vspace{0.5cm}
L.B. Ostrovsky, M.E. Zhukovskii~\footnote{Moscow Institute of Physics and Technology, Laboratory of Advanced Combinatorics and Network Applications}
\end{center}

\vspace{1cm}

\begin{center}
Abstract
\end{center}

%\vspace{0.3cm}

We study asymptotical probabilities of first order and monadic second order properties of Erdos--Renyi random graph $G(n,n^{-a})$. The random graph obeys FO (MSO) zero-one $k$-law if for any first order (monadic second order) formulae it is true for $G(n,n^{-a})$ with probability tending to $0$ or to $1$. Zero-one $k$-laws are well studied only for the first order language and $a<1$. We obtain new zero-one $k$-laws (both for first order and monadic second order languages) when $a>1$. Proofs of these results are based on the existed study of first order equivalence classes and our study of monadic second order equivalence classes. The respective results are of interest by themselves.

\vspace{1cm}

\section{Logic of the random graph}
\label{zero_one_intro}

In the paper, we study asymptotical probabilities of first order and monadic second order properties of Erd\H{o}s--R\'{e}nyi random graph $G(n,p)$~\cite{Bollobas,Janson,Strange,Survey}. Recall that edges in this graph on the set of vertices $V_n=\{1,\ldots,n\}$ appear independently with probability $p$ (i.e., for any undirected graph $H=(V_n,E)$ without loops and multiple edges the equality ${\sf P}(G(n,p)=H)=p^{|E|}(1-p)^{{n\choose 2}-|E|}$ holds).

Formulae in the first order language of graphs (first order formulae)~\cite{Logic1,Logic2,Survey,Veresh,Strange} are constructed using relational symbols $\sim$ (the symbol of adjacency) and $=$; logical connectivities $\neg,\Rightarrow,\Leftrightarrow,\vee,\wedge$; variables $x,y,x_1, \ldots$ (that express vertices of a graph); and quantifiers $\forall,\exists$. Monadic second order formulae~\cite{Muller,Tysk} are built of the above symbols of the first order language and variables $X,Y,X_1,\ldots$ that express unary predicates. Following \cite{Logic1,Logic2,Survey,Veresh}, we call a number of nested quantifiers in the longest chain of nested quantifiers of a formula $\phi$ {\it the quantifier depth} $q(\phi)$. For example, the formula
$$
 (\forall X \,\,([\exists x_1\exists x_2\,\,
 (X(x_1)\wedge(\neg(X(x_2))))]\Rightarrow
 [\exists y \exists z\,\,(X(y)\wedge(\neg(X(z)))\wedge(y\sim z))]))
$$
has quantifier depth $3$ and expresses the property of being connected. It is known that this property is not expressed by a first order formula (see, e.g.,~\cite{Survey}).

We say that $G(n,p)$ {\it obeys FO zero-one law} ({\it MSO zero-one law}) if for any first order formula (monadic second order formula) it is either true asymptotically almost surely (a.a.s.) or false a.a.s. (as $n\to\infty$). In 1988, S.~Shelah and J.~Spencer~\cite{Shelah} proved the following zero-one law for the random graph $G(n,n^{-\alpha})$.
\begin{theorem}
Let $\alpha>0$. The random graph $G(n,n^{-\alpha})$ does not obey FO zero-one law if and only if either $\alpha\in(0,1]\cap\mathbb{Q}$ or $\alpha=1+1/l$ for some integer $l$.
\label{FO01}
\end{theorem}
Obviously, there is no MSO zero-one law when even FO zero-one law does not hold. In 1993, J.~Tyszkiewicz~\cite{Tysk} proved that $G(n,n^{-\alpha})$ does not obey MSO zero-one law for irrational $\alpha\in(0,1)$ also. When $\alpha>1$ and does not equal to any of $1+1/l$ MSO zero-one law holds. The last statement simply follows from standart arguments from the theory of logical equivalence. We believe, this result is known. Unfortunately, we did not find it in the related papers. So, we give the proof in Section~\ref{proof1}. Below, we state the general result on MSO zero-one law for $G(n,n^{-\alpha})$.
\begin{theorem}
Let $\alpha>0$. The random graph $G(n,n^{-\alpha})$ does not obey MSO zero-one law if and only if either $\alpha\in(0,1]$ or $\alpha=1+1/l$ for some integer $l$.
\label{MSO01}
\end{theorem}

For a formula $\phi$ consider the set $S(\phi)$ of $\alpha$ such that $G(n,n^{-\alpha})$ does not obey the zero-one law for the fixed formula $\phi$. Both theorems do not give any explanation of how the set $S(\phi)$ depends on $\phi$ (or even on a quantifier depth of this formula). However, better insight into an asymptotical behavior of probabilities of the properties expressed by first order and monadic second order formulae is given by zero-one $k$-laws (see Section~\ref{k-laws}), which are well studied only for the first order language and $\alpha\leq 1$ (see, e.g.,~\cite{Survey}). In the presented paper, we obtain new zero-one $k$-laws (both for first order and monadic second order languages) when $\alpha>1$ and give their statements in Section~\ref{k-laws}. Proofs of these results are based on the existed study of first order equivalence classes and our study of monadic second order equivalence classes (see Section~\ref{elementary}). The respective results are of interest by themselves.

\section{Logical equivalence}
\label{elementary}

%We write $G\models\phi$ if the formula $\phi$ is true on $G$.

For two graphs $G$ and $H$ and any positive integer $k$, the notation $G\equiv_k^{\FO;\,\gr}H$ denotes that any first order formula $\phi$ with $q(\phi)\leq k$ is true on both $G$ and $H$ or false on both $G$ and $H$. The notation $G\equiv_k^{\MSO;\,\gr} H$ is defined similarly. Obviously, $\equiv_k^{\FO;\,\gr}$ and $\equiv_k^{\MSO;\,\gr}$ are both equivalence relations on the set of all graphs. Moreover, for every $k$ there are only finitely many equivalence classes (see, e.g.,~\cite{Logic2}) and an upper bound for the cardinality $r^{\FO;\,\gr}_k$ of the set of all $\equiv_k^{\FO;\,\gr}$-equivalence classes $\mathcal{R}^{\FO;\,\gr}_k$ is known~\cite{Verbitsky} and given below. Let $T(s)$ be the tower function: $T(s)=2^{T(s-1)}$, $T(1)=2$. Let $\log^*(k)=\min\{i:\,T(i)\geq k\}$.
\begin{theorem}
$$
 r^{\FO;\,\gr}_k\leq T(k+2+\log^*(k))+O(1).
$$
\label{FO_equiv}
\end{theorem}
Similarly, $r^{\MSO;\,\gr}_k$ and $\mathcal{R}^{\MSO;\,\gr}_k$ are defined. In this paper, we prove a similar result for $\equiv^{\MSO;\,\gr}_k$-equivalence (the proof is given in Section~\ref{proof2}).
\begin{theorem}
For any positive integer $k$,
$$
 r^{\MSO;\,\gr}_k\leq T(k+2+\log^*(k)).
$$
\label{MSO_equiv}
\end{theorem}
Note that this result is stronger than Theorem~\ref{FO_equiv}.

In order to prove the results on zero-one laws from Section~\ref{k-laws}, we also need an extension of the above theory to the case of rooted trees. Recall that {\it a rooted tree} $T_R$ is a tree with one distinguished vertex $R$, which is called {\it the root}. If $R,\ldots,x,y$ is a simple path in $T_R$, then $x$ is called a parent of $y$ and $y$ is called a child of $x$. The first order language for rooted trees has a constant symbol $R$ (for the root) and the parent--child relation $P(x,y)$. For two rooted trees $T_R$ and $T'_{R'}$ and any positive integer $k$, the notations $T_R\equiv_k^{\FO;\,\tr}T'_{R'}$, $T_R\equiv_k^{\FO;\,\tr}T'_{R'}$, $r^{\FO;\,\tr}_k$, $r^{\MSO;\,\tr}_k$, $\mathcal{R}^{\FO;\,\tr}_k$, $\mathcal{R}^{\MSO;\,\tr}_k$ are defined in the same way as for graphs. The following result is proven in~\cite{Verbitsky}.
\begin{theorem}
$$
 r^{\FO;\,\tr}_k\leq T(k+2+\log^*(k))+O(1).
$$
For any $A\in \mathcal{R}^{\FO;\,\tr}_k$, the following inequality holds:
$$
 \min_{T_R\in A} |V(T_R)|\leq T(k+4+\log^*(k))+O(1).
$$
\label{FO_tr_equiv}
\end{theorem}
In this paper, we prove a similar result for $\equiv^{\MSO;\,\tr}_k$-equivalence (the proof is given in Section~\ref{proof3}).
\begin{theorem}
Let $k\geq 4$ be an integer. Then
\begin{equation}
 r^{\MSO;\,\tr}_k\leq T(k+2+\log^*(k)).
\label{trees_number}
\end{equation}
For any $A\in \mathcal{R}^{\MSO;\,\tr}_k$, the following inequality holds:
$$
 \min_{T_R\in A} |V(T_R)|\leq T(k+3+\log^*(k+1)).
$$
\label{MSO_tr_equiv}
\end{theorem}
Note that this result is stronger than Theorem~\ref{FO_tr_equiv}.

\section{Zero-one $k$-laws}
\label{k-laws}

By Theorem~\ref{FO01}, for any rational $\alpha\in(0,1]$ there is a first order formula $\phi$ which is true on $G(n,p)$ with probability which asymptotics either does not exist or does not equal to $0$ or $1$. Obviously, this statement is not true for formulae with bounded quantifier depth. We say that $G(n,p)$ {\it obeys FO zero-one $k$-law} ({\it MSO zero-one $k$-law}) if for any first order formula (monadic second order formula) with quantifier depth at most $k$ it is either true a.a.s. or false a.a.s. (as $n\to\infty$). In~\cite{Zhuk_01_ext}--\cite{Zhuk_law}, the following zero-one $k$-laws are proven.
\begin{theorem}
For any $k\geq 3$ and any $\alpha\in(0,1/(k-2))$ the random graph $G(n,n^{-\alpha})$ obeys FO zero-one $k$-law. If $\alpha=1/(k-2)$, then $G(n,n^{-\alpha})$ does not obey FO zero-one $k$-law. For any $k\geq 4$ and any $\alpha\in(1-1/(2^k-2),1)$ the random graph $G(n,n^{-\alpha})$ obeys FO zero-one $k$-law. If $\alpha=1-1/(2^k-2)$, then $G(n,n^{-\alpha})$ does not obey FO zero-one $k$-law.
\label{old_k-law}
\end{theorem}

In this paper, we consider the very sparse case $\alpha>1$. From Theorems~\ref{FO01},~\ref{MSO01}, $G(n,n^{-\alpha})$ obeys both zero-one $k$-laws if $\alpha\neq 1+1/l$ for any positive integer $l$. In Sections~\ref{proof4},~\ref{proof5} we prove the following result.
\begin{theorem}
Let $l$ be positive integer, $\alpha=1+1/l$.
\begin{itemize}
\item[$\cdot$]Let $k \ge 4$ be an arbitrary integer. If $l \ge T(k + \log^*(k+1) + 3)$, then the random graph $G(n, n^{-\alpha})$ obeys MSO zero-one $k$-law.
\item[$\cdot$]Let $k \ge 7$ be an arbitrary integer. If $l\leq 2T(k-4)$, then the random graph $G(n,n^{-\alpha})$ does not obey FO zero-one $k$-law.
\end{itemize}
\label{new_k-law}
\end{theorem}

\section{Proofs}

We start from notations and auxiliary statements (Section~\ref{pre}). Remind that the simple case ($\alpha>1$) of the statement of Theorem~\ref{MSO01} is known. We give its proof in Section~\ref{proof1}, because we did not find it in the related papers. Theorem~\ref{new_k-law} is proven in Sections~\ref{proof4},~\ref{proof5}. These proofs are based on Theorems~\ref{MSO_equiv},~\ref{MSO_tr_equiv}. Proofs of the latter results can be found in Sections~\ref{proof2},~\ref{proof3} respectively.

\subsection{Preliminaries}
\label{pre}

Throughout this paper, if $G$ is a graph then we denote its vertex set by $V(G)$ and its edge set by $E(G)$ (i.e. $G=(V(G),E(G))$). The {\it distance} between vertices $u$ and $v$ in a connected graph is the minimum length of a path connecting $u$ and $v$. It is denoted by $d(u,v)$. For a disconnected graph, the distance between vertices in different components equals $\infty$. {\it The eccentricity} of a vertex $v$ is $e(v) = \max_{u\in V(G)} d(v, u)$. {\it The diameter} and {\it the radius} of $G$ are $d(G)=\max_{v\in V(G)} e(v)$ and $r(G) = \min_{v\in V(G)}e(v)$ respectively. A vertex $v$ is {\it central}, if $e(v)=r(G)$.

For a rooted tree $T_R$, we call the eccentricity of its root {\it the depth of $T_R$}. The relation of being {\it a descendant} is the transitive and reflexive closure of the relation of being a child. If $v\in V(T_R)$, then $T_R(v)$ denotes the subtree of $T_R$ spanned by the set of all descendants of $v$ and rooted at $v$.

For a graph $G$ and a formula $\phi$, we write $G\models\phi$ if $\phi$ is true on $G$.\\

In this section, we review well-known statements (and prove new, see Section~\ref{pre_logic}, Lemma~\ref{Verb_Lem} and Lemma~\ref{Verb_Lem2}) from the random graph theory (Section~\ref{pre_graphs}) and the model theory (Section~\ref{pre_logic}), which are exploited in our proofs.

\subsubsection{Small subgraphs of the random graph}
\label{pre_graphs}

Consider a graph $G$ on $v$ vertices and $e$ edges. Denote the number of automorphisms of $G$ by $a(G)$. The fraction $\rho(G) = \frac{e}{v}$ is called {\it the density} of $G$. The graph $G$ is called {\it strictly balanced} if for any proper subgraph $H \subset G$ the inequality $\rho(H) < \rho(G)$ holds. Let $N_G$ be the number of copies of $G$ in $G(n,p)$. Let $G$ be strictly balanced. In~\cite{Erdos}, a threshold probability for the property of containing $G$ was obtained.

\begin{theorem}
If $p\gg n^{-1/\rho(G)}$, then a.a.s. in $G(n,p)$ there is a copy of $G$. Moreover, $\frac{N_G}{{\sf E}N_G}\stackrel{{\sf P}}\rightarrow 1$. If $p\ll n^{-1/\rho(G)}$, then a.a.s. in $G(n,p)$ there are no copies of $G$.
\label{threshold}
\end{theorem}

Moreover, in~\cite{Bol_small} an asymptotical distribution of $N_G$ in the threshold was found.
\begin{theorem}
If $p=cn^{-1/\rho(G)}$, then $N_G\stackrel{d}\rightarrow\Pois(c^{e}/a(G))$.
\label{poisson}
\end{theorem}

A threshold probability for the property of being connected is stated in the following result (see, e.g.,~\cite{Bollobas,Janson}).
\begin{theorem}
Let $c<1$. If for $n$ large enough $p(n)<c\frac{\ln n}{n}$, then a.a.s. $G(n,p)$ is not connected. If for $n$ large enough $p(n)>(1/c)\frac{\ln n}{n}$, then a.a.s. $G(n,p)$ is connected.
\label{connectedness}
\end{theorem}
From these results it follows that for any positive integer $l$ and any $1+1/(l+1)<\alpha<1+1/l$ the following three properties hold.
\begin{itemize}
\item[T1] The random graph $G(n,n^{-\alpha})$ is a forest a.a.s.
\item[T2] A.a.s. any component of $G(n,n^{-\alpha})$ has at most $l+1$ vertices.
\item[T3] For any integer $K$, a.a.s. for any tree $T$ on at most $l+1$ vertices there are at least $K$ components in $G(n,n^{-\alpha})$ which are isomorphic to $T$.
\end{itemize}
Moreover, for any positive integer $l$ and $\alpha=1+1/l$ the properties T1 and T2 hold. Moreover,
\begin{itemize}
\item[T4] For any tree $T$ on $l+1$ vertices the probability of containing $T$ tends to $1-e^{-1/a(T)}$.
\end{itemize}

\subsubsection{Ehrenfeucht game}
\label{pre_logic}

The main tool of all the above results is Ehrenfeucht game~\cite{Logic1,Logic2,Ehren}, \cite{Muller}--\cite{Strange}, \cite{Veresh,Survey}. We start from the general first order theory on arbitrary finite structures. Consider the first order language consisting of arbitrary relational symbols $P_1,\ldots,P_m$ of arities $a_1,\ldots,a_m$ respectively and constant symbols $R_1,\ldots,R_s$. The game $\EHR^{\FO}(A,B,k)$ is played on structures $A$ (with distinguished elements $R^A_1,\ldots,R^A_s$) and $B$ (with distinguished elements $R^B_1,\ldots,R^B_s$) of the above vocabulary. There are two players (Spoiler and Duplicator) and a fixed number of rounds $k$. At the $\nu\mbox{-}$th round ($1 \leq \nu \leq k$) Spoiler chooses either an element $x_{\nu}$ of $A$ or an element $y_{\nu}$ of $B$. Duplicator chooses an element of the other structure. In the case of monadic second order logic, players can choose subsets as well. Similarly, at the $\nu\mbox{-}$th round ($1 \leq \nu \leq k$) of the game $\EHR^{\MSO}(A,B,k)$ Spoiler chooses any structure of $A$ and $B$. Say, he chooses $A$. Then he either chooses an element $x_{\nu}$ or a subset $X_{\nu}$ of $A$. If an element is chosen, Duplicator chooses an element of $B$. Otherwise, Duplicator chooses a subset of $B$.

In $\EHR^{\FO}$, at the end of the game the elements $x_{1},...,x_{k}$ of $A$, $y_{1},...,y_{k}$ of $B$ are chosen. Denote $x_{k+1}=R^A_1,y_{k+1}=R^B_1,\ldots,x_{k+s}=R^A_s,y_{k+s}=R^B_s$. Duplicator wins if and only if the following property holds.

\begin{enumerate}

\item[$\cdot$] For any $i\in\{1,\ldots,m\}$ and $\nu_1,\ldots,\nu_a\in\{1,\ldots,k+s\}$ (where $a=a_i$), $P_i(x_{\nu_1},\ldots,x_{\nu_a})\Leftrightarrow P_i(y_{\nu_1},\ldots,y_{\nu_a})$.

\end{enumerate}

In $\EHR^{\MSO}$, at the end of the game elements $x_{i_1},...,x_{i_t}$ of $A$, $y_{i_1},...,y_{i_t}$ of $B$ and subsets $X_{j_1},...,X_{j_{k-t}}$ of $A$, $Y_{j_1},...,Y_{j_{k-t}}$ of $B$ are chosen. Duplicator wins if and only if the following two properties hold.

\begin{enumerate}

\item[$\cdot$] For any $i\in\{1,\ldots,m\}$ and $\nu_1,\ldots,\nu_a\in\{i_1,\ldots,i_t,k+1,\ldots,k+s\}$ (where $a=a_i$), $P_i(x_{\nu_1},\ldots,x_{\nu_a})\Leftrightarrow P_i(y_{\nu_1},\ldots,y_{\nu_a})$.

\item[$\cdot$] For any $\nu\in\{i_1,\ldots,i_t,k+1,\ldots,k+s\}$ and $\mu\in\{j_1,\ldots,j_{k-t}\}$,
$x_{\nu}\in X_{\mu}\Leftrightarrow y_{\nu}\in Y_{\mu}$.

\end{enumerate}

Further, we follow definitions and notations from~\cite{Verbitsky} (see Section 3). Here, we reformulate these definitions for the monadic second order theory on arbitrary structures. Let $t,l$ be nonnegative integers, $t+l\leq k$. Consider two structures $A$ and $B$, elements $u_1,\ldots,u_t$ of $A$, $v_1,\ldots,v_t$ of $B$ and subsets $U_1,\ldots,U_l$ of $A$, $V_1,\ldots,V_l$ of $B$. Let elements $u_1,\ldots,u_t,v_1,\ldots,v_t$ be chosen in the first $t$ rounds of $\EHR^{\MSO}(A,B,k)$, sets $U_1,\ldots,U_l,V_1,\ldots,V_l$ be chosen in the next $l$ rounds. Suppose Duplicator wins $\EHR^{\MSO}(A,B,t+l)$. Denote $\overline{u}=(u_1,\ldots,u_t)$, $\overline{v}=(v_1,\ldots,v_t)$, $\overline{U}=(U_1,\ldots,U_l)$, $\overline{V}=(V_1,\ldots,V_l)$. We write $(A,\overline{u},\overline{U})\equiv_k^{\MSO}(B,\overline{v},\overline{V})$ if Duplicator has a winning strategy in the rest $k-t-l$ rounds. Obviously, this $\equiv_k^{\MSO}$ is an equivalence relation (the $\equiv_k^{\FO}$-equivalence is defined in the same way, see~\cite{Verbitsky}). The {\it $k$-Ehrenfeucht value} of $(A,\overline{u},\overline{U})$ is the $\equiv_k^{\MSO}$-equivalence class it belongs to. We let $\EHR(k,t,l)$ denote the set of all $k$-Ehrenfeucht values for structures of the considered vocabulary with marked $s$ elements and $m$ subsets. Set $\EHR(k)=\EHR(k,0,0)$.

Our proofs of Theorems~\ref{MSO_equiv},~\ref{MSO_tr_equiv} (see Sections~\ref{proof2},~\ref{proof3} respectively) are based on an extension of Lemma 3.2 from~\cite{Verbitsky} to the monadic second order language, which is stated below.

\begin{lemma}
Consider two cases.
\begin{itemize}
\item[$\cdot$] If $t+l=k$, then
\begin{equation}
|\EHR(k,t,l)|\leq 2^{\sum_{i=1}^m a_i!{{s+t}\choose {a_i}}} 2^{(s+t)l}.
\label{last}
\end{equation}
\item[$\cdot$] If $t+l<k$, then
\begin{equation}
|\EHR(k,t,l)|\leq 2^{|\EHR(k,t+1,l)|+|\EHR(k,t,l+1)|}.
\label{recur}
\end{equation}
\end{itemize}
\label{Verb_Lem}
\end{lemma}
{\it Proof.} Let $t + l = k$. Then the $\equiv_k^{\MSO}$-class of $A$ with marked $\overline{u},\overline{U}$ is determined by memberships of $u_1,\ldots,u_t,R_1,\ldots,R_s$ in $U_1,\ldots,U_l$, and the induced substructure $A|_{\{u_1,\ldots,u_t,R_1,\ldots,R_s\}}$. Therefore, we get~(\ref{last}).

If $s+m<k$, then
$(A,\overline{u},\overline{U})\equiv_k^{\MSO}(B,\overline{v},\overline{V})$ if and only if
\begin{equation}
\forall u\in A\exists v \in B\quad (A, (u_1,\ldots,u_m,u), \overline{U}) \equiv_{k}^{\MSO} (B, (v_1,\ldots,v_m,v), \overline{V}),
\label{F1}
\end{equation}
\begin{equation}
\forall U \subset A\exists V\subset B\quad
(A, \overline{u}, (U_1,\ldots,U_s,U)) \equiv_{k}^{\MSO} (B, \overline{v}, (V_1,\ldots,V_s,V)),
\label{F2}
\end{equation}
and vice versa. Obviously~(\ref{F1}),~(\ref{F2}) hold if and only if
\begin{itemize}
\item[$\cdot$] the set of $k$-Ehrenfeucht values of the structure $A$ with marked elements $u_1,\ldots,u_t,u$ and subsets $U_1,\ldots,U_l$ over all $u\in A$ and the set of $k$-Ehrenfeucht values of the structure $B$ with marked elements $v_1,\ldots,v_t,v$ and subsets $V_1,\ldots,V_l$ over all $v\in B$ coincide;
\item[$\cdot$] the set of $k$-Ehrenfeucht values of the structure $A$ with marked elements $u_1,\ldots,u_t$ and subsets $U_1,\ldots,U_l,U$ over all $U\subset A$ and the set of $k$-Ehrenfeucht values of the structure $B$ with marked elements $v_1,\ldots,v_t$ and subsets $V_1,\ldots,V_l,V$ over all $V\subset B$ coincide.
\end{itemize}
Therefore, the $k$-Ehrenfeucht value of the structure $A$ with marked elements $u_1,\ldots,u_t$ and subsets $U_1,\ldots,U_l$ is defined by the set of $k$-Ehrenfeucht values of the structure $G$ with marked elements $u_1,\ldots,u_t,u$ and subsets $U_1,\ldots,U_l$ over all $u\in A$ and the set of $k$-Ehrenfeucht values of the graph $G$ with marked elements $u_1,\ldots,u_t$ and subsets $U_1,\ldots,U_l,U$ over all $U\subset A$. This leads to Equation~\ref{recur}. Lemma is proven.\\

For graphs (vocabulary consists of two relational symbols $\sim,=$ and no constant symbols, see Section~\ref{zero_one_intro}) and rooted trees (two relational symbols $P,=$ and one constant symbol $R$ are considered, see Section~\ref{elementary}), the bound~(\ref{last}) can be strengthened.

\begin{lemma}
Let $t+l=k$. For graphs, $|\EHR(k,t,l)|\leq 2^{k^2-k}$. For rooted trees and $k\geq 5$, $\log_2 |\EHR(k,t,l)|\leq 2^k-2$. For rooted trees and $k=4$, $|\EHR(4,t,l)|\leq 3\cdot 2^{13}$.
\label{Verb_Lem2}
\end{lemma}
{\it Proof.} As relations $\sim,=$ are symmetric, for graphs,
$$
 |\EHR(k,t,l)|\leq 2^{2{t\choose 2}+tl}=2^{t(2k-t-1)}\leq 2^{k(k-1)}.
$$
For rooted trees, we consider two cases: vertices $u_1,\ldots,u_t,R$ are either pairwise distinct or not. Note that there are at most $(t+1)^{t-1}3^t$ directed forests on $t+1$ labeled vertices (there are $(t+1)^{t-1}$ trees, each edge can be either removed or directed in one of two ways). Therefore, in the first case (for distinct $u_1,\ldots,u_t,R$), the number of $k$-Ehrenfeucht values of $(A,\overline{u},\overline{U})$ is at most
$$
 (t+1)^{t-1}3^t 2^{(t+1)l}=2^{(t-1)\log_2(t+1)+t\log_2 3+(t+1)l}\leq 2^{g_1(k)},
$$
where $g_1(k)=\left(\frac{k+1}{2}\right)^2+(k-1)\log_2(k+1)+k\log_2 3$. As the number of ways of determining the equality relation on $u_1,\ldots,u_t,R$ is at most $2^{{{t+1}\choose 2}}$, in the second case ($u_1,\ldots,u_t,R$ are not pairwise distinct), the number of $k$-Ehrenfeucht values is less than
$$
 2^{{{t+1}\choose 2}}t^{t-2}3^{t-1}2^{tl}\leq 2^{{{t+1}\choose 2}+(t-2)\log_2 t+(t-1)\log_2 3+tl}\leq 2^{g_2(k)},
$$
where $g_2(k)=\frac{(k+1)k}{2}+(k-2)\log_2 k+(k-1)\log_2 3$. Therefore, for $k\geq 5$, we get the following bound:
$$
 |\EHR(k,t,l)|\leq 2^{g_1(k)}+2^{g_2(k)}<2^{2^k-2}.
$$
The last inequality holds because, first, for $k=5$, both $g_1(k),g_2(k)$ are less than $2^k-3$, and, second, $2^k-3$ growths faster than $g_1(k)$ and $g_2(k)$ for $k>5$.

Let $k=4$. In the first case (for distinct $u_1,\ldots,u_t,R$), the number of $k$-Ehrenfeucht values of $(A,\overline{u},\overline{U})$ is at most
$$
 \max_{t\in\{0,1,2,3,4\}}(t+1)^{t-1}3^t 2^{(t+1)(4-t)}=3^4 5^3.
$$
There are exactly $F(t)=1,3,19,201$ directed forests on $t=1,2,3,4$ labeled vertices respectively. Moreover, for $t=1,2,3,4$, the number of ways of determining the equality relation on $u_1,\ldots,u_t,R$ in such a way that these vertices are not pairwise distinct equals $E(t)=1,4,14,51$ respectively. Finally, if $u_1,\ldots,u_t,R$ are not pairwise distinct, then, for $t=1,2,3,4$, the number of ways of determining memberships of $u_1,\ldots,u_t,R$ in $U_1,\ldots,U_{4-t}$ is at most $M(t)=8,16,8,1$ respectively. Therefore, in the second case ($u_1,\ldots,u_t,R$ are not pairwise distinct), the number of $k$-Ehrenfeucht values is at most
$$
 \max_{t\in\{0,1,2,3,4\}}F(t)E(t)M(t)=201\cdot 51.
$$
Finally,
$$
 |\EHR(4,t,l)|\leq 3^4 5^3+201\cdot 51=20376<3\cdot 2^{13}.
$$
Lemma is proven.\\

For graphs and vocabulary consisting of two relational symbols $\sim,=$ and no constant symbols, we denote the games $\EHR^{\FO}(G,H,k)$, $\EHR^{\MSO}(G,H,k)$ on graphs $G,H$ by
$$
\EHR^{\FO;\,\gr}(G,H,k),\quad\EHR^{\MSO;\,\gr}(G,H,k)
$$
respectively. In 1960~\cite{Ehren}, A.~Ehrenfeucht proved that there is the following connection between Ehrenfeucht games and logical equivalences (we give the statement confirming to the case of graphs).
\begin{theorem}
For any positive integer $k$ and any graphs $G,H$,
$$
 G\equiv^{\FO;\,\gr}_k H\Leftrightarrow\text{Duplicator has a winning strategy in }\EHR^{FO;\,\gr}(G,H,k),
$$
$$
 G\equiv^{\MSO;\,\gr}_k H\Leftrightarrow\text{Duplicator has a winning strategy in }\EHR^{MSO;\,\gr}(G,H,k).
$$
\label{ehren_original}
\end{theorem}
We need the following well-known corollary from this theorem (see, e.g.,~\cite{Muller,Janson,Strange,Survey}).
\begin{corol}
Let $k$ be any positive integer. The random graph $G(n,p)$ obeys FO zero-one $k$-law if and only if asymptotically almost surely Duplicator has a winning strategy in $\EHR^{\FO;\,\gr}(G(n,p),G(m,p),k)$ as $n,m\to\infty$. The random graph $G(n,p)$ obeys MSO zero-one $k$-law if and only if asymptotically almost surely Duplicator has a winning strategy in $\EHR^{\MSO;\,\gr}(G(n,p),G(m,p),k)$ as $n,m\to\infty$.
\label{ehren}
\end{corol}

The following two standart facts about $\equiv^{\MSO;\,\gr}_k$-equivalence are essential tools in our arguments. They follow from Theorem~\ref{ehren_original} (see the proofs in~\cite{Compton}, Theorems 2.2 and 2.3). Denote by $A\sqcup B$ a disjoint union of graphs $A,B$: $A\sqcup B=(V(A)\sqcup V(B),E(A)\sqcup E(B))$.

\begin{lemma}
Let $k$ be a positive integer. If $H_1\equiv^{\MSO;\,\gr}_k H_2$, $G_1\equiv^{\MSO;\,\gr}_k G_2$, then $H_1\sqcup G_1\equiv^{\MSO;\,\gr}_k H_2\sqcup G_2$. The same holds for $\equiv^{\FO;\,\gr}_k$-equivalence.
\label{L1}
\end{lemma}

Denote by $aG$ a disjoint union of $a$ copies of a graph $G$.

\begin{lemma}
For any positive integer $k$ there is a positive integer $a$ such that for any integer $b\geq a$ and any graph $G$ we have $bG\equiv^{\MSO;\,\gr}_k aG$.
\label{L2}
\end{lemma}

In the proofs, we also exploit Ehrenfeucht game on rooted trees. We denote the games $\EHR^{\FO}(T_R,T'_{R'},k)$, $\EHR^{\MSO}(T_R,T'_{R'},k)$ on rooted trees $T_R,T'_{R'}$ by $\EHR^{\FO;\,\tr}(T_R,T'_{R'},k)$, $\EHR^{\MSO;\,\tr}(T_R,T'_{R'},k)$ respectively (two relational symbols $P,=$ and one constant symbol $R$ are considered). Another particular case of Ehrenfeucht Theorem is stated below.
\begin{theorem}
For any positive integer $k$ and any rooted trees $T_R,T'_{R'}$,
$$
 T_R\equiv^{\FO;\,\tr}_k T'_{R'}\Leftrightarrow\text{Duplicator has a winning strategy in }\EHR^{FO;\,\tr}(T_R,T'_{R'},k),
$$
$$
 T_R\equiv^{\MSO;\,\tr}_k T'_{R'}\Leftrightarrow\text{Duplicator has a winning strategy in }\EHR^{MSO;\,\tr}(T_R,T'_{R'},k).
$$
\end{theorem}

\subsection{Proof of Theorem~\ref{MSO01}: $\alpha>1$}
\label{proof1}

Let $l$ be a positive integer.

We start from $\alpha\in(1+1/(l+1),1+1/l)$ and $p=n^{-\alpha}$.  Let $k$ be a positive integer. Consider $a=a(k)$ from Lemma~\ref{L2}. Let $\mathcal{T}$ be the set of all pairwise nonisomorphic trees on at most $l+1$ vertices. Consider the forest
$$
T_0=\bigsqcup_{T\in\mathcal{T}}aT.
$$
From Properties T1, T2, T3 it follows that there exist $K_{T}(n),T\in\mathcal{T}$, such that a.a.s.
$$
G(n,p)=\bigsqcup_{T\in\mathcal{T}}K_T(n)T
$$
and $K_T(n)\to\infty$ as $n\to\infty$ for any $T\in\mathcal{T}$. From Lemma~\ref{L2}, for $n$ large enough $K_T(n)T\equiv_k^{\MSO;\,\gr}aT$ for any $T\in\mathcal{T}$. Finally, from Lemma~\ref{L1}, for such $n$
$$
T_0\equiv_k^{\MSO;\,\gr}\bigsqcup_{T\in\mathcal{T}}K_T(n)T.
$$
Therefore, a.a.s. $G(n,p)\equiv_k^{\MSO;\,\gr}T_0$. This means that a.a.s. for any monadic second order $\phi$ with $q(\phi)=k$ a.a.s. $G(n,p)\models\phi$ if and only if $T_0\models\phi.$ Thus, $G(n,p)$ obeys MSO zero-one $k$-law. As $k$ is arbitrary, MSO zero-one law holds.

If $\alpha=1+1/l$ and $p=n^{-\alpha}$, then $G(n,p)$ does not obey FO zero-one law. Therefore, MSO zero-one law fails as well.

If $\alpha>2$, then a.a.s. in $G(n,n^{-\alpha})$ there are no edges (i.e. $G(n,n^{-\alpha})$ is the union of isolated vertices). From Lemma~\ref{L2}, in this case MSO zero-one law holds.

\subsection{Proof of Theorem~\ref{MSO_equiv}}
\label{proof2}

Fix a positive integer $k$. To avoid trivialities, we assume $k\geq 2$. Set $\tilde T(1)=2^{k^2-k}$, $\tilde T(i)=2^{2\tilde T(i-1)}$, $\hat T(1)=2^{2^k}$, $\hat T(i)=2^{\hat T(i-1)}$. As $2^k\geq k^2$ and for a positive $x$, $1+2^x<2^x+2^x=2^{x+1}$, we get the inequalities $\tilde T(i)\leq\hat T(i)$, $i\in\{1,\ldots,k+1\}$. Moreover, $\hat T(1)\leq 2^{2^{T(\log^*(k))}}=T(2+\log^*(k))$. For any nonnegative integer $\beta$ such that $\beta\leq k$ set $f(k,\beta)=\max_{t\in\{0,1,\ldots,\beta\}}|\EHR(k,t,\beta-t)|$. By Lemma~\ref{Verb_Lem} and Lemma~\ref{Verb_Lem2},
$$
 |\EHR(k)|\leq 2^{2f(k,1)}\leq 2^{2\cdot 2^{2f(k,2)}}\leq\ldots\leq\tilde T(k+1)\leq T(k+2+\log^*(k)).
$$
Finally, by Theorem~\ref{ehren_original}, we get $r_k^{\MSO;\,\gr}=|\EHR(k)|\leq T(k+2+\log^*(k))$.

\subsection{Proof of Theorem~\ref{MSO_tr_equiv}}
\label{proof3}

Fix an integer $k\geq 5$. Set $\tilde T(1)=2^{2^k-2}$, $\tilde T(i)=2^{2\tilde T(i-1)}$, $\hat T(1)=k-1+2^{2^k-1}$, $\hat T(i)=2^{\hat T(i-1)}$. As in the previous section, for any $i\in\{1,\ldots,k+1\}$, we have $\tilde T(i)\leq\hat T(i)$. Moreover, $\hat T(1)<2^{2^k}\leq 2^{2^{T(\log^*{k})}}=T(2+\log^*(k))$. Therefore, by Lemma~\ref{Verb_Lem} and Lemma~\ref{Verb_Lem2},
$$
 |\EHR(k)|\leq \tilde T(k+1)\leq \hat T(k+1)<T(k+2+\log^*(k)).
$$
Let $k=4$. By Lemma~\ref{Verb_Lem} and Lemma~\ref{Verb_Lem2},
$$
 \log_2\log_2|\EHR(k)|\leq 1+2^{1+2^{1+6\cdot 2^{13}}}<2^{2^{3+6\cdot 2^{13}}}<2^{2^{2^{16}}}=T(6)=T(k+\log^*(k)).
$$
By Theorem~\ref{ehren_original}, for any $k\geq 4$ we get $r_k^{\MSO;\,\tr}=|\EHR(k)|\leq T(k+2+\log^*(k))$.\\

The second statement of Theorem~\ref{MSO_tr_equiv} follows the statement about representatives of $\equiv_k^{\MSO;\,\tr}$-classes which is an extension of Lemma 8.6 and Lemma 8.7 of~\cite{Verbitsky} to the monadic second order language and stated below. Set $z=2^kf(k,1)\ldots f(k,k)$ (see the notation in Section~\ref{proof2}).

\begin{lemma}
Let $k$ be a positive integer. Fix any $A\in\mathcal{R}_k^{\MSO}$. Let $T_R\in A$ be a rooted tree with a minimal order over all rooted trees in $A$. Then
\begin{itemize}
\item each vertex of $T_R$ has at most $zf(k,0)$ children;
\item the depth of $T_R$ is at most $f(k,0)$.
\end{itemize}
\label{represent}
\end{lemma}
{\it Proof.} For $\equiv^{\MSO}_k$-equivalence class, we say that is has {\it type $m$}, if any its representative contains $t+l=m$ vertices and subsets (see Section~\ref{pre_logic}). For all $m\in\{0,1,\ldots,k\}$, set $z(k,m)=1$. Moreover, for all $i\in\{1,\ldots,k\}$, $m\in\{0,1,\ldots,i\}$, set
\begin{equation}
 z(i-1,m)=z(i,m)+z(i,m+1)f(k,m+1).
\label{introducing_z}
\end{equation}
Obviously, for all $i\in\{0,1,\ldots,k-1\}$, $m\in\{0,1,\ldots,i\}$,
$$
 z(i,m)\leq 2^{k-i}f(k,m+1)\ldots f(k,k).
$$
Therefore, $z(0,0)\leq z$.

Consider an arbitrary vertex $y$ of $T_R$. Let $w_1,\ldots,w_s$ be children of $y$. Let $k$-Ehrenfeucht values of $T_R(w_1),\ldots, T_R(w_s)$ equal $a_1,\ldots,a_s$ respectively. Suppose that for some $i\in\{1,\ldots,s\}$ (say, $i=1$) more than $z$ values of $a_1,\ldots,a_s$ equal $a_i$. Let $w_1,\ldots,w_{z'}$ be all children of $y$ such that $k$-Ehrenfeucht values of $T_R(w_1),\ldots, T_R(w_{z'})$ equal $a_1$. Consider the rooted tree $T^-_R$ which is obtained from $T_R$ by removing all but $z$ subtrees rooted at children of $y$ such that their $k$-Ehrenfeucht values equal $a_1$ (say, the subtrees $T_R(w_{z+1}),\ldots,T_R(w_{z'})$).

Prove that $T_R\equiv^{\MSO;\,\tr}_k T^-_R$. By Theorem~\ref{ehren_original}, the trees are equivalent if and only if Duplicator has a winning strategy in $\EHR^{\MSO;\,\tr}(T_R,T^-_R,k)$. We do not consider choices of vertices and subsets of
$$
V(T_R)\setminus V(T_R(w_{1})\cup\ldots\cup T_R(w_{z'}))
$$
by Spoiler, because for such choices the strategy of Duplicator is trivial. Set
$$
 \mathcal{C}_{0}(0)=\{T_R(w_1),\ldots,T_R(w_{z'})\},\quad
 \mathcal{C}^{-}_0(0)=\{T^-_R(w_1),\ldots,T^-_R(w_{z})\}.
$$
Consider the first round of $\EHR^{\MSO;\,\tr}(T_R,T^-_R,k)$. If Spoiler chooses a vertex $u\neq y$ (say, in a tree $T_R(w_i)$), then Duplicator chooses any vertex $v$ such that $(T_R(w_i),u)\equiv_k^{\MSO}(T^-_R(w_j),v)$ (say, in a tree $T^-_R(w_j)$). He can do this, because all trees from $\mathcal{C}^{-}_0$ and $\mathcal{C}_0$ have the same $k$-Ehrenfeucht value. Set
$$
 \mathcal{C}_{1;1}(1)=\{(T_R(w_i),u)\},\quad\mathcal{C}_{1;1}(0)=\{T_R(w_1),\ldots,T_R(w_{z'})\}\setminus \{T_R(w_i)\},
$$
$$
 \mathcal{C}^-_{1;1}(1)=\{(T^-_R(w_j),v)\},\quad
 \mathcal{C}^-_{1;1}(0)=\{T^-_R(w_1),\ldots,T^-_R(w_{z})\}\setminus \{T^-_R(w_j)\},
$$
$$
 \gamma_1(0)=1, \gamma_1(1)=1.
$$
Finally, if Spoiler chooses a subset $U$ (say, $U=U_1\cup\ldots\cup U_{z'}$, $U_i$ is (maybe, empty) subset of $V(T_R(w_i))$), then $\mathcal{C}_{1;\cdot}(1)$ are $\equiv_k^{\MSO}$-equivalence classes on pairs $(T_R(w_i),U_i)$ such that $U_i$ are not empty, and $\mathcal{C}_{1;\cdot}(0)$ is the $\equiv_k^{\MSO}$-equivalence class (maybe, empty) on $T_R(w_i)$ such that $U_i$ are empty. Obviously, the number of equivalence classes $\mathcal{C}_{1;\cdot}(1)$ equals $\gamma_1(1)\leq f(k,1)$, the number of equivalence classes $\mathcal{C}_{1;\cdot}(0)$ equals $\gamma_1(0)\leq 1$. Duplicator constructs a subset $V$ in the following way. For each $m\in\{0,1\}$ and each $j$ from $1$ to $\gamma_1(m)$, if $|\mathcal{C}_{1;j}(m)|\leq z(1,m)$, then Duplicator takes next $|\mathcal{C}_{1;j}(m)|$ trees of $T^-_R(w_1),\ldots,T^-_R(w_{z})$ and chooses their subsets (empty for $m=0$) such that each tree with the respective subset is in the same $\equiv_k^{\MSO}$-equivalence class with any representative of $\mathcal{C}_{1;j}(m)$. If $|\mathcal{C}_{1;j}(m)|> z(1,m)$, then Duplicator takes next at least $z(1,m)$ trees and their subsets in the same way. As
$$
z'>z\geq z(0,0)=z(1,0)+z(1,1)f(k,1),
$$
the subset $V$ exists.

Consider the $i$-th round, $i\in\{2,\ldots,k\}$. Let, after the previous round, nonnegative integers $\gamma_{i-1}(0),\ldots,\gamma_{i-1}(i-1)$ and $\equiv^{\MSO}_k$-equivalence classes $\mathcal{C}_{i-1;j}(m)$, $\mathcal{C}^-_{i-1;j}(m)$, $m\in\{0,\ldots,i-1\}$, $j\in\{1,\ldots,\gamma_{i-1}(j)\}$, be chosen (each class consists of rooted trees from $\mathcal{C}_0$ and $\mathcal{C}^-_0$ respectively with chosen $m$ vertices and subsets). %Let
%$$
%\gamma_{i-1}(j)\leq 2^{j-1}f(k,1)\ldots f(k,j)
%$$
%for $j\in\{1,\ldots,i-1\}$ and $\gamma_{i-1}(0)\leq 1$. Moreover,
For any $m\in\{0,\ldots,i-1\}$ and any $j\in\{1,\ldots,\gamma_{i-1}(m)\}$, all elements of $\mathcal{C}^{-}_{i-1;j}(m)$ and $\mathcal{C}_{i-1;j}(m)$ have the same $k$-Ehrenfeucht value. Finally, for any $m\in\{0,\ldots,i-1\}$ and any $j\in\{1,\ldots,\gamma_{i-1}(m)\}$,
\begin{equation}
\text{either }|\mathcal{C}^{-}_{i-1;j}(m)|=|\mathcal{C}_{i-1;j}(m)|,\text{ or }\min\{|\mathcal{C}^{-}_{i-1;j}(m)|,|\mathcal{C}_{i-1;j}(m)|\}\geq z(i-1,m).
\label{lem_recur_main}
\end{equation}
If Spoiler chooses a vertex $u\neq y$ (say, in a tree $T_R(w_i)$ with chosen $\hat m$ vertices $\overline{u}$ and subsets $\overline{U}$, $(T_R(w_i),\overline{u},\overline{U})\in\mathcal{C}_{i-1;1}(\hat m)$), then Duplicator chooses any vertex $v$ such that $(T_R(w_i),\overline{u},u,\overline{U})\equiv_k^{\MSO}(T^-_R(w_{\tilde i}),\overline{v},v,\overline{V})$ (where $(T^-_R(w_{\tilde i}),\overline{v},\overline{V})$ is from $\mathcal{C}^-_{i-1;1}(\hat m)$). Set
$$
 \mathcal{C}_{i;1}(\hat m)=\mathcal{C}_{i-1;1}(\hat m)\setminus\{(T_R(w_i),\overline{u},\overline{U})\},\quad
 \mathcal{C}_{i;j}(\hat m)=\mathcal{C}_{i-1;j}(\hat m),\quad j\in\{2,\ldots,\gamma_{i-1}(\hat m)\},
$$
$$
 \mathcal{C}_{i;1}(\hat m+1)=\{(T_R(w_i),\overline{u},u,\overline{U})\},\quad
 \mathcal{C}_{i;j+1}(\hat m+1)=\mathcal{C}_{i-1;j}(\hat m+1),\quad j\in\{1,\ldots,\gamma_{i-1}(\hat m+1)\},
$$
$$
 \mathcal{C}_{i;j}(m)=\mathcal{C}_{i-1;j}(m),\quad m\in\{0,\ldots,i-1\}\setminus\{\hat m,\hat m+1\},\quad j\in\{2,\ldots,\gamma_{i-1}(m)\}.
$$
In the same way, the classes $\mathcal{C}^-_{i;j}(m)$ are defined. Moreover,
$\gamma_i(m)=\gamma_{i-1}(m)$ for all $m\in\{0,\ldots,i\}\setminus\{\hat m+1\}$ and $\gamma_i(\hat m+1)=\gamma_{i-1}(\hat m+1)+1$ (here, we assume $\gamma_{i-1}(i)=0$). Finally, if Spoiler chooses a subset $U$ (say, $U=U_1\cup\ldots\cup U_{z'}$, $U_l$ is (maybe, empty) subset of $V(T_R(w_l))$), then classes $\mathcal{C}_{i;\cdot}(m+1)$ are obtained by dividing each $\mathcal{C}_{i-1;j}(m)$ into at most $f(k,m+1)$ $\equiv^{\MSO}_k$-equivalence classes with respect to nonempty subsets $U_l$. Elements of a class $\mathcal{C}_{i-1;j}(m)$ with empty subsets $U_l$ form (at most one) class of the same type $m$. %Therefore, the number of equivalence classes $\mathcal{C}_{i;\cdot}$ of type $j$ equals
%$$
% \gamma_i(j)\leq\gamma_{i-1}(j-1)f(k,j)+\gamma_{i-1}(j)\leq 2^{j-1}f(k,1)\ldots f(k,j)
%$$
%for $j\in\{1,\ldots,i-1\}$, $\gamma_i(i)\leq 2^{i-2}f(k,1)\ldots f(k,i)$ and $\gamma_i(0)\leq 1$. Obviously, %$$
%\gamma_i=\sum_j\gamma_i(j)\leq.
%$$
Duplicator constructs a subset $V$ in the following way. For each $m$ from $0$ to $i$ and each $j$ from $1$ to $\gamma_i(m)$, he finds the class $\mathcal{C}_{i-1;\tilde j}(\hat m)\supset\mathcal{C}_{i;j}(m)$. If $|\mathcal{C}_{i;j}(m)|\leq z(i,m)$, then Duplicator takes next $|\mathcal{C}_{i;j}(m)|$ elements of the class $\mathcal{C}^-_{i-1;\tilde j}(\hat m)$ and chooses (maybe, empty) subsets of the respective rooted trees such that each tree with chosen subsets (maybe, accounting for the new one which is nonempty) and vertices is in the same $\equiv_k^{\MSO}$-equivalence class with any representative of $\mathcal{C}_{i;j}(m)$. If $|\mathcal{C}_{i;j}(m)|> z(i,m)$, then Duplicator takes next at least $z(i,m)$ elements of the class $\mathcal{C}^-_{i-1;\tilde j}(\hat m)$ in the same way. The subset $V$ exists, because (\ref{introducing_z}) with $m=\hat m$ and (\ref{lem_recur_main}) hold.

Obviously, this strategy of Spoiler is winning. Therefore, $T^-_R\in A$ and $|V(T^-_R)|<|V(T_R)|$, that contradicts with the assumptions. As the number of different classes among $a_1,\ldots,a_s$ is at most $f(k,0)$, the first statement of Lemma~\ref{represent} is proven.

Suppose there is a vertex $y\in V(T_R)$ such that $d(R,y)>f(k,0)$. Consider the path $R,\ldots,y$. By the definition of $f(k,0)$, in this path there are distinct nonroot vertices $w_1,w_2$ such that $T_R(w_1)\equiv^{\MSO}_k T_R(w_2)$. Assume $d(R,w_1)<d(R,w_2)$. Consider a rooted at $R$ tree $T_R^{-}$ which is obtained from $T_R$ by the result of replacing $T_R(w_1)$ by $T_R(w_2)$ (the root $w_1$ is replacing by the root $w_2$). It is easy to see that $T_R\equiv^{\MSO}_k T_R^{-}$. By Theorem~\ref{ehren_original}, $T_R\equiv^{\MSO;\,\tr}_k T^-_R$. Therefore, $T^-_R\in A$ and $|V(T^-_R)|<|V(T_R)|$, that contradicts with the assumptions. Lemma is proven.\\

Fix $A\in\mathcal{R}_k^{\MSO}$. From Lemma~\ref{represent} and the first statement of Theorem~\ref{MSO_tr_equiv}, we get
$$
 \min_{T_R\in A}|V(T_R)|\leq \sum_{i=0}^{f(k,0)}{zf(k,0)}^i<(z f(k,0))^{f(k,0)+1}\leq T(k+3+\log^*(k+1)).
$$
Theorem is proven.

\subsection{Proof of Theorem~\ref{new_k-law}: MSO $k$-law}
\label{proof4}

Let $k\geq 4$, $l\geq T(k+\log^* (k+1)+3)$, $\alpha=1+\frac{1}{l}$, $p=n^{-\alpha}$. Consider $a=a(k)$ from Lemma~\ref{L2}. Let $\mathcal{T}$ be the set of all pairwise nonisomorphic trees on at most $l$ vertices and $\mathcal{T}^+$ be the set of all pairwise nonisomorphic trees on exactly $l+1$ vertices. Consider the forest
$$
T_0=\bigsqcup_{T\in\mathcal{T}}aT.
$$
From Properties T1, T2, T3 it follows that there exist $K_{T}(n),T\in\mathcal{T}\cup\mathcal{T}^+$, such that a.a.s.
$$
G(n,p)=\bigsqcup_{T\in\mathcal{T}\cup\mathcal{T}^+}K_T(n)T
$$
and $K_T(n)\to\infty$ as $n\to\infty$ for any $T\in\mathcal{T}$. Obviously, for any $\equiv_k^{\MSO;\,\tr}$-equivalence class, trees in this class are in one $\equiv_k^{\MSO;\,\gr}$-equivalence class. Therefore, from Theorem~\ref{MSO_tr_equiv}, for any $\hat T\in\mathcal{T}^+$ there exists $T(\hat T)\in\mathcal{T}$ such that $T(\hat T)\equiv_k^{\MSO;\,\gr}\hat T$. Denote $\mathcal{T}(\mathcal{T}^+)$ the set of all $T\in\mathcal{T}$ such that there exists $\hat T$ with $T=T(\hat T)$. For any $T\in\mathcal{T}(\mathcal{T}^+)$, denote $\mathcal{T}^+(T)$ the set of trees $\hat T$ from $\mathcal{T}^+$ such that $T=T(\hat T)$. From Lemma~\ref{L2}, for $n$ large enough $\tilde K_T(n) T\equiv_k^{\MSO;\,\gr}aT$ for any $T\in\mathcal{T}$, where $\tilde K_T(n)=K_T(n)+\sum_{\hat T\in\mathcal{T}^+(T)}K_{\hat T}(n)$ for all $T\in\mathcal{T}(\mathcal{T}^+)$ and $\tilde K_T(n)=K_T(n)$ for all the others $T\in\mathcal{T}$. Finally, from Lemma~\ref{L1}, for such $n$
$$
 T_0\equiv_k^{\MSO;\,\gr}\bigsqcup_{T\in\mathcal{T}}\tilde K_T(n)T\equiv_k^{\MSO;\,\gr}\left(\bigsqcup_{T\in\mathcal{T}}K_T(n)T\right)\sqcup
 \left(\bigsqcup_{{\hat T}\in\mathcal{T}^+}K_{\hat T}(n){\hat T}\right).
$$
Therefore, a.a.s. $G(n,p)\equiv_k^{\MSO;\,\gr}T_0$. This means that a.a.s. for any monadic second order $\phi$ with $q(\phi)=k$ a.a.s. $G(n,p)\models\phi$ if and only if $T_0\models\phi.$ Thus, $G(n,p)$ obeys MSO zero-one $k$-law.

\subsection{Proof of Theorem~\ref{new_k-law}: FO $k$-law}
\label{proof5}

In the proof, we follow definitions and notations from~\cite{Verbitsky}, Section 6. We give them below. If $w$ is a child of $u  \in V(T_v)$, then $T_{v}(w)$ is a {\it $u$-branch} of $T_v$. A rooted tree $T_v$ is called {\it diverging} if for any its vertex $u$ all $u$-branches of $T_v$ are pairwise nonisomorphic. A tree is called {\it diverging} if for any its central vertex $v$ the rooted tree $T_v$ is diverging. We say that rooted trees $T_v$ and $S_u$ are {\it isomorphic} if there is an isomorphism $T\rightarrow U$ which maps $v$ to $u$.

We use the following result, which is proved in~\cite{Verbitsky} (see Lemma 6.13).

\begin{lemma}
Let $i \ge 3$. For every $n$ such that $2 i + 2 \le n \le 2 T(i - 1)+1$ there exists a diverging tree of order $n$ and radius $i + 1$.
\label{construct}
\end{lemma}

For each $l\in\{1,\ldots,8\}$, it is easy to construct a first order formula with quantifier depth at most $7$ which is true on $G(n,n^{-1-1/l})$ with asymptotical probability in $(0,1)$ (using the properties T1 and T4). Let $l\geq 7$.

Fix $k\geq 7$. We start from $2T(k-4)\geq l\geq 2k-5$. Set $p=n^{-\alpha}$, $\alpha=1+1/l$. By Lemma~\ref{construct}, there exists a diverging tree $S$ of order $l + 1$ and radius $k - 2$. By T1, T2 and T4, with some asymptotical probability $c\in(0,1)$ in $G(n,p)$ there is a component isomorphic to $S$. Moreover, with asymptotical probability $1-c$ in $G(n,p)$ there is no copy of $S$.

Consider two forests $A$ and $B$ such that there is a component $S^A$ in $A$ isomorphic to $S$, and there is no copy of $S$ in $B$. Let us prove that Spoiler has a winning strategy in $\EHR^{\FO;\,\gr}(A, B, k)$. In the first round, Spoiler chooses a central vertex $x_1$ in $S^A$. Duplicator chooses a vertex $y_1$ in a component $S^B$ of $B$.

Suppose that $d(S^A)< d(S^B)$ (if $d(S^A)>d(S^B)$, Spoiler applies the same strategy). In the next two rounds, Spoiler chooses $y_2$ and $y_3$ (not necessary different from $y$) in $S^B$ such that $d(y_2,y_3) = d(S^A) + 1$ and $d(y_1,y_2) \le d(S^A)$. Duplicator chooses vertices $x_2,x_3$. If $x_2\notin V(S^A)$, then there is a winning strategy of Spoiler in next $\lceil\log_2(d(S^A))\rceil\leq\lceil\log_2(2r(S^A))\rceil=1+\log_2(k-2)\leq k-3$ rounds (see,~e.g.,~\cite{Verbitsky},~Lemma 6.4). If $x_2\in V(S_A)$, then $d_{S_A}(x_2, x_3) \ne d(y_2,y_3)=d(S^A)+1$. Therefore,  there is a winning strategy of Spoiler in next $\lceil\log_2(d(S^A)+1)\rceil\leq\lceil\log_2(2k-3)\rceil\leq k-3$ rounds.

Let $d(S^A) = d(S^B)$. Further, we apply the strategy of Spoiler from Lemma 6.7 and Lemma 6.8 in~\cite{Verbitsky}. Note that this strategy is winning in a game on two trees. However, in the main part of this strategy, Spoiler in each round chooses a vertex which is adjacent to one of the vertices chosen in the previous rounds. Therefore, Duplicator can not change a tree (i.e., this strategy is also winning in a game on forests). To make our proof self-contained, we sketch this strategy in the paragraph below.

If $y_1$ is not a central vertex of $S^B$, then as in the previous cases Spoiler has a winning strategy in next at most $1+\lceil r(S^A)\rceil<k-1$ rounds. Let $y_1$ be a central vertex of $S^B$. Suppose that the tree $S^B$ is diverging. Let us prove that for any $i\in\{1,\ldots,k\}$ either Spoiler wins in the $i$-th round or chosen vertices $x_1, x_2, \ldots, x_i$ and $y_1, y_2, \ldots, y_i$ form simple paths, and the rooted trees  $S^A_{x_1}(x_i)$ and $S^B_{y_1}(y_i)$ are not isomorphic. For $i=1$, this is already proven. Assume that for some $i\in\{1,\ldots,k\}$ this statement is also proven. Consider the round $i+1$. If only one vertex of $x_i,y_i$ is a leaf, Spoiler wins. Let $x_i,y_i$ be not leafs. As the rooted trees $S^A_{x_1}$ and $S^B_{y_1}$ are diverging, $x_i$-branches of $S^A_{x_1}$ are pairwise nonisomorphic ($y_i$-branches of $S^B_{y_1}$ are pairwise nonisomorphic as well), see~\cite{Verbitsky}, Lemma 6.2. Without loss of generality, there is a vertex $x_{i+1}$ which is a child of $x_i$ such that there is no $y_i$-branch isomorphic to $S^A_{x_{i+1}}$. Spoiler chooses $x_{i+1}$, and the statement is proven. Therefore, Spoiler wins in at most $r(S^A)+1=k-1$ rounds. If the tree $S^B$ is not diverging, consider a vertex $t$ of $S^B$ such that $S^B_{y_1}(t)$ is not a diverging tree but any its $t$-branch is diverging. Let two isomorphic $t$-branches be rooted at $z_1,z_2$. Spoiler selects the path $y_1,y_2,\ldots,t,z_1$, Duplicator's response is  $x_1,x_2,\ldots,x_i$. If depths of $S^B_{y_1}(z_1)$ and $S^A_{x_1}(x_i)$ are distinct, Spoiler prolongs a path corresponding to a smaller depth and wins in at most $r(S^A) + 1\leq k-1$ rounds. If depths are equal and $S^B_{y_1}(z_1),S^A_{x_1}(x_i)$ are not isomorphic, then the winning strategy of Spoiler in $r(S^A)+1$ rounds is given above. Finally, if $S^B_{y_1}(z_1),S^A_{x_1}(x_i)$ are isomorphic, in the round $i+1$ Spoiler chooses $z_2$ and Duplicator responses in $x_{i+1}$, a child of $x_{i-1}$. As $S^A_{x_1}(x_{i-1})$ is a diverging tree, the trees $S^A_{x_1}(x_{i+1})$ and $S^B_{y_1}(z_2)$ are not isomorphic. Analogously, Spoiler wins in at most $r(S^A) + 2\leq k$ rounds.

From the above,
$$
 \mathrm{liminf}_{n\rightarrow\infty}{\sf P}(\text{Spoiler has a winning  strategy in }\EHR(G(n,p(n)),G(m,p(m)),k))\geq c(1-c).
$$
By Theorem~\ref{ehren}, MSO zero-one $k$-law fails when $2k-5\leq  l\leq 2T(k-4)$ and $k\geq 7$.  Therefore, MSO zero-one $k$-law fails for all $l\in\{9,\ldots,2k-6\}$ as well. Theorem is proven.

\section{Acknowledgments}

The authors are grateful to Professor Oleg Verbitsky for helpful discussions on
the Workshop on Logic and Random Graphs in the Lorentz Center (August 31 --- September
4, 2015).

\end{document}